\documentclass[sn-mathphys,Numbered]{sn-jnl}

\usepackage{graphicx}%
\usepackage{multirow}%
\usepackage{amsmath,amssymb,amsfonts}%
\usepackage{amsthm}%
\usepackage{mathrsfs}%
\usepackage[title]{appendix}%
\usepackage{xcolor}%
\usepackage{textcomp}%
\usepackage{manyfoot}%
\usepackage{booktabs}%
\usepackage{algorithm}%
\usepackage{algorithmicx}%
\usepackage{algpseudocode}%
\usepackage{listings}%
\usepackage{enumerate}
\usepackage{setspace}

\theoremstyle{thmstyleone}%
\newtheorem{theorem}{Theorem}[section]%

\theoremstyle{thmstyletwo}%
\newtheorem{remark}{Remark}[section]%

\theoremstyle{thmstylethree}%
\newtheorem{definition}{Definition}[section]%

\theoremstyle{thmstylefour}%
\numberwithin{equation}{section}
\begin{document}
	
	\title[Finslerian exponentially harmonic Liouvielle theorem]{A Liouville-type theorem for Finslerian exponentially harmonic functions}
	
	
	\author*[1]{\fnm{Bin} \sur{Shen}}\email{shenbin@seu.edu.cn}
	
	
	
	\affil*[1]{\orgdiv{School of Mathematics}, \orgname{Southeast University}, \orgaddress{\street{Dongnandaxue Road 2}, \city{Nanjing}, \postcode{211189}, \state{Jiangsu}, \country{China}}}


	
	\abstract{ In this manuscript, we investigate the exponentially harmonic equation on noncompact forward complete Finsler metric measure spaces. We demonstrate that this Finslerian equation represents a critical point of an exponential energy functional. Furthermore, we establish that any bounded solution to this equation is constant, provided that the mixed weighted Ricci curvature is nonnegative and certain additional non-Riemannian tensors are bounded.}

	
	
	
	\maketitle

	\section{Introduction}
	Exponentially harmonic maps between Riemannian manifolds were initially introduced by Eells and Lemaire in 1990 \cite{EL1990}. Subsequently, Duc and Eells investigated the existence and regularity properties of exponentially harmonic functions \cite{DE1991}. Hong established the equivalence between harmonic and exponentially harmonic maps on closed Riemannian manifolds under conformal metric transformations \cite{Hong1992-1}. Furthermore, he demonstrated the Liouville theorem for exponentially harmonic functions defined on complete Riemannian manifolds \cite{Hong1992-2}. However, it should be noted that in general cases, there exists no inclusion relationship between the spaces of harmonic functions and exponentially harmonic functions \cite{HY1993}. The extension of exponentially harmonic maps to Finsler manifolds was subsequently explored by Chiang in \cite{Chiang2018}.
	
	In this manuscript, we investigate Finslerian exponentially harmonic functions. Analogous to the Riemannian case, a function $u$ defined on a Finsler metric measure space $(M,F,\mu)$ is said to be exponentially harmonic if it satisfies 
	\begin{eqnarray}\label{equ-expharmonic}
		\tilde\Delta_{\mu}u:=\mathrm{div}_{\mu}(V(u)Du)=0,
	\end{eqnarray}
	where $V(u)=\exp(\frac12F^{*2}(Du))$ represents the energy density of $u$. Such a function arises as a critical point of the exponential energy functional (see Theorem \ref{thm-EL}). Finsler metric measure spaces encompass a richer set of geometric tensors compared to Riemannian manifolds (cf. \cite{ShenLec}). It is important to emphasize that \eqref{equ-expharmonic} is not merely a quasilinear elliptic equation, as the underlying space is an anisotropic and asymmetric manifold. Within this framework, we establish the following theorem.
	\begin{theorem}
		Let $(M,F,\mu)$ be a forward complete $n$-dimensional Finsler metric measure space with finite misalignment $\alpha$.  Assume that the mixed weighted Ricci curvature $^mRic^{\infty}$ of $M$ is nonnegative, and that the $S$-curvature as well as the non-Riemannian curvatures $U$, $\mathcal{T}$ and $divC=FC^{i}_{\,jk\vert i}dx^j\otimes dx^k$ satisfy the norm bound
		$F^{-1}|S|+F(U)+F(\mathcal{T})+\|divC\|_{HS(V)} \leq K_0,$
		for any $(x,V)\in SM$.
		Then, any bounded exponentially harmonic function $u$ on $M$ is constant.
	\end{theorem}
	
	\section{Related concepts and notations of the Finsler metric measure spaces}
	
	A Finsler metric measure space is a triple $(M,F,\mu)$, where $M$ is a differentiable manifold equipped with a Finsler metric $F$ and a Borel measure $\mu$.
	The \emph{Cartan tensor} is defined by
	$$C(X,Y,Z):=C_{ijk}X^iY^jZ^k=\frac{1}{4}\frac{\partial^3F^2(x,y)}{\partial y^i\partial y^j\partial y^k}X^iY^jZ^k,$$
	for any local vector fields $X,Y,Z$. 
	We always adopt the Chern connection, which is uniquely determined by
	\begin{align*}
		\nabla_XY-\nabla_YX&=[X,Y];\\
		Z(g_y(X,Y))-g_y(\nabla_ZX,Y)-&g_y(X,\nabla_ZY)=2C_y(\nabla_Zy,X,Y),
	\end{align*}
	for any $X,Y,Z\in TM\setminus\{0\}$, where $C_y$ is the Cartan tensor. The coefficients of the Chern connection are locally expressed as $\Gamma^i_{jk}(x,y)$ in natural coordinates. These coefficients induce the \emph{spray coefficients} as $G^i=\frac12\Gamma^i_{jk}y^jy^k$. The spray is given by 
	\begin{eqnarray}
		G=y^i\frac{\delta}{\delta x^i}=y^i\frac{\partial}{\partial x^i}-2G^i\frac{\partial}{\partial y^i},
	\end{eqnarray}
	where $\frac{\delta}{\delta x^i}=\frac{\partial}{\partial x^i}-N^j_i\frac{\partial}{\partial y^j}$, and the nonlinear connection coefficients $N^i_j$ are locally derived from the spray coefficients as $N^i_j=\frac{\partial G^i}{\partial y^j}$.
	By convention, the horizontal Chern derivative is denoted by $``|"$ and the vertical Chern derivative by $``;"$. Let $\hat D$ denote the Levi-Civita connection of the induced Riemannian metric $\hat g=g_Y$, and let $\{e_i\}$ be an orthonormal basis with respect to metric $g(x,Y)$ at a point $x$, where $Y$ is a fixed reference vector field. Furthermore, let $\{E_i\}$ be the local vector fields obtained by parallel transporting $\{e_i\}$ in a neighborhood of $x$ on $M$. In \cite{Shen2024}, the tensor $U$ is defined as 
	\begin{eqnarray}
		U_y(W)=g(x,W)(U(y,W),W),
	\end{eqnarray} 
	for any local vector field $W$, with 
	\begin{eqnarray}\label{def-UyW}
		U(y,W)=\sum_{i=1}^n (D^W_{e_i}E_i-\hat D_{e_i}E_i)
	\end{eqnarray}
	representing a vector field on the sphere bundle, where $Y$ denotes the geodesic extension field of $y$ in a neighborhood of $x$. 
	The Riemannian structure is not invariant under changes in the reference vector. In \cite{Shen2024}, an important constant on a Finsler manifold, called the misalignment, is introduced. It is also provided in \cite{W.Zhao:2013} as a universe constant $\Lambda_F$.
	\begin{definition}[\cite{Shen2024}]
		The misalignment of a Finsler manifold $(M,F)$ at a point $x$ is defined by
		\begin{eqnarray}
			\alpha(x)=\sup_{V,W,Y\in S_xM}\frac{g_V(Y,Y)}{g_W(Y,Y)}.
		\end{eqnarray}
		Furthermore, the global misalignment of the Finsler metric is defined by
		\begin{eqnarray}\label{def-alpha}
			\alpha=\sup_{x\in M}\alpha_M(x)=\sup_{x\in M}\sup_{V,W,Y\in S_xM}\frac{g_V(Y,Y)}{g_W(Y,Y)}.
		\end{eqnarray}
	\end{definition}
	In \cite{Shen2024}, several characterizations of the misalignment are provided. Specifically, a Finsler manifold $(M,F)$ is a Riemannian manifold if and only if $\alpha_M=1$. Additionally, a Finsler manifold $(M,F)$ exhibits uniform convexity and uniform smoothness if and only if it satisfies finite misalignment.

	A \emph{Finsler metric measure space} $(M,F,d\mu)$ is a Finsler manifold equipped with a given measure $\mu$. 
	In local coordinates $\{x^i\}_{i=1}^n$, the volume form can be expressed as $d\mu=\sigma(x)dx^1\wedge\cdots\wedge dx^n$, where $\sigma(x)\geq 0$. For any $y\in T_xM\setminus\{0\}$, define 
	$$\tau(x,y):=\log\frac{\sqrt{\det g_{ij}(x,y)}}{\sigma(x)},$$
	which is called the \emph{distortion} of $(M,F,\mu)$. The definition of the \emph{S-curvature} is as follows.
	\begin{definition}[\cite{ShenLec}]\label{def-S}
		Let $(M,F,d\mu)$ be a Finsler metric measure space. For any point $x\in M$, let $\gamma=\gamma(t)$ be a forward geodesic starting from $x$ with the initial tangent vector $\dot\gamma(0)=y$. The S-curvature of $(M,F,d\mu)$ is defined as 
		\begin{eqnarray*}
			S(x,y):=\frac{d}{dt}\tau(\gamma(t),\dot\gamma(t))|_{t=0}=\frac{d}{dt}\left(\frac12\log\det(g_{ij})-\log\sigma(x)\right)(\gamma(t),\dot\gamma(t))|_{t=0}.
		\end{eqnarray*}
	\end{definition}
	Definition \ref{def-S} implies that the S-curvature measures the rate of change of the distortion along the geodesic in the direction of $y$. Inspired by the definition of the $T$-curvature in \cite{ShenLec}, \cite{Shen2024} introduced the difference of $\nabla\tau$ on the tangent sphere, denoted by $\mathcal{T}$. That is,
	\begin{eqnarray}\label{def-TVW}
		\mathcal{T}(Y,W):=\nabla^V\tau(Y)-\nabla^W\tau(W),
	\end{eqnarray}
	for vector fields $Y,W$ on $M$. Locally, it can be expressed as $\mathcal{T}(Y,W)=\mathcal{T}_i(Y,W)d x^i$, where 
	\begin{eqnarray}
		\mathcal{T}_i=\frac{\delta}{\delta x^i}\tau_k(Y)-\frac{\delta}{\delta x^i}\tau_k(W).
	\end{eqnarray}
	Clearly, $\mathcal{T}(Y,W)$ is anti-symmetric about $Y$ and $W$, i.e., $\mathcal{T}(Y,W)=-\mathcal{T}(W,Y)$ for any nonvanishing $Y$ and $W$.
	
	Using the S-curvature, we can now define the\emph{ weighted Ricci curvature} as follows.	
	\begin{definition}[\cite{Ohta2014}]
		Given a unit vector $Y\in T_xM$ and a positive number $k$, the weighted Ricci curvature is defined by
		\begin{eqnarray}
			Ric^k(Y):=\begin{cases}
				&Ric(x,Y)+\dot{S}(x,Y)\quad\mbox{if } S(x,Y)=0 \mbox{ and } k=n \mbox{ or if }k=\infty;\\
				&-\infty\quad\quad\quad\quad\quad\quad\quad\quad\quad\quad\quad\quad\mbox{if }S(x,Y)\neq0 \mbox{ and if } k=n;\\
				&Ric(x,Y)+\dot{S}(x,Y)-\frac{S^2(x,Y)}{k-n}\quad\quad\quad\quad\quad\quad\quad\mbox{ if }n<k<\infty,
			\end{cases}
		\end{eqnarray}
		where the derivative is taken along the geodesic starting from $x$ in the direction of $Y$.
	\end{definition}
	
	Moreover, the \emph{mixed weighted Ricci curvature} has been defined in \cite{Shen2024}, denoted by $^mRic^k$.
	\begin{definition}[\cite{Shen2024}]\label{Def-mwrc}
		Given two unit vectors $Y,W\in T_xM$ and a positive number $k$, the mixed weighted Ricci curvature $^mRic^k(Y,W)=\,^mRic^k_W(Y)$ is defined by
		\begin{eqnarray}
			^mRic^k_{W}(Y):=\begin{cases}
				&\mathrm{tr}_WR_{Y}(Y)+\dot{S}(x,Y)\quad\mbox{if } S(x,Y)=0 \mbox{ and } k=n \mbox{ or if }k=\infty;\\
				&-\infty\quad\quad\quad\quad\quad\quad\quad\quad\quad\quad\quad\quad\mbox{if }S(x,Y)\neq0 \mbox{ and if } k=n;\\
				&\mathrm{tr}_WR_{Y}(Y)+\dot{S}(x,Y)-\frac{S^2(x,Y)}{k-n}\quad\quad\quad\quad\quad\quad\quad\mbox{ if }n<k<\infty,
			\end{cases}
		\end{eqnarray}
		where the derivative is taken along the geodesic starting from $x$ in the direction of $Y$, and $\mathrm{tr}_WR_{Y}(Y)=g^{ij}(W)g_{Y}(R_Y(e_i,Y)Y,e_j)$ represents the trace of the flag curvature with respect to $g(x,W)$.
	\end{definition}
	
	\begin{remark}
		The weighted Ricci curvature is a special case of the mixed weighted Ricci curvature, i.e., $Ric^k(Y)=\,^mRic^k_Y(Y)$. Specifically, when $W=0$, the mixed weighted Ricci curvature reduces to the weighted Ricci curvature, i.e., $\,^mRic^k_W(Y)=Ric^k(Y)$.
	\end{remark}
	
	Defining the function $\mathfrak{ct}_c(r)$ as
	\begin{eqnarray}
		\mathfrak{ct}_c(r)=\begin{cases}
			\sqrt{c}\cot\sqrt{c}r,\quad\,\,\,\quad c>0,\\
			1/r, \quad\quad\quad\quad\quad\,\,\,\,  c=0,\\
			\sqrt{-c}\coth\sqrt{-c}r, \quad c<0.
		\end{cases}
	\end{eqnarray}
	The following significant Laplacian comparison theorem was first obtained in \cite{Shen2024}, and has since been successfully applied to the study of nonlinear Finslerian parabolic equations \cite{ShenZhu2025}.
	\begin{theorem}[\cite{Shen2024}]\label{thm-LapComp-1}
		Let $(M,F,\mu,\alpha)$ be an $n$-dimensional forward complete Finsler metric measure space with finite misalignment $\alpha$.
		Let $r$ denote the forward distance function and $V$ a nontrivial vector field on $M$. Suppose that for some $N\in (n,\infty)$  the mixed weighted Ricci curvature satisfies $^mRic^N_{\nabla r}\geq-K$ with $K\geq0$, and the non-Riemannian curvatures $U$, $\mathcal{T}$ and $divC=FC^{i}_{\,\,\,\,jk\vert i}dx^\otimes dx^k$
		satisfy the norm bounds by $F(U) + F^*(\mathcal{T} ) + \|divC\|_{HS(V)} \leq K_0$. Then, by setting $l=K/C(N,\alpha)$, where $C(N,\alpha)=N+(\alpha-1)n-\alpha$, the nonlinear Laplacian of $r$ with reference vector $V$ satisfies 
		\begin{eqnarray*}
			\Delta^Vr\leq C(N,\alpha)\mathfrak{ct}_{-l}(r)+C_0,
		\end{eqnarray*}
		wherever $r$ is $C^2$, where $C_0$ is a constant depending on $\alpha$ 
		and the bound of non-Riemannian curvatures $K_0$.
	\end{theorem}
	\begin{remark}
		According to Definition \ref{Def-mwrc}, the curvature condition $^mRic^N_{\nabla r}\geq-K$ in Theorem \ref{thm-LapComp-1} can be replaced by $^mRic^{\infty}_{\nabla r}\geq-K$ and $F^{-1}|S|\leq K_0$.
	\end{remark}
	We consistently omit the reference vector in the gradient and the Laplacian when the reference direction aligns with the gradient of the function under consideration.

	\section{Finslerian exponentially harmonic functions}
	
	
	Suppose $(M,F,d\mu)$ is a forward complete Finsler metric measure space. Let $u$ be a function defined on $M$. The exponential energy is given by the following functional.
	\begin{eqnarray}
		\mathbb{E}(u)=\int_M\exp\frac{F^{*2}(Du)}{2}d\mu.
	\end{eqnarray}
	
	\begin{theorem}\label{thm-EL}
		The Euler-Lagrange operator associated with $\mathbb{E}$ is the quasilinear strictly elliptic operator 
		\begin{eqnarray}\label{def-FEH}
			\tilde\Delta_{\mu}u:=\mathrm{div}_{\mu}(\mathcal{V}(u)Du)=e^{-\sigma}\frac{\delta}{\delta x^i}(g^{ij}(\nabla u)e^{\sigma+\frac12F^{*2}(Du)}u_j),
		\end{eqnarray}
		where $\mathcal{V}(u)=\exp(\frac12F^{*2}(Du))$ is the energy density of $u$.
	\end{theorem}
	\begin{proof}
		Let $v$ be a local function with compact support. A direct variational calculation shows that
		\begin{eqnarray}
			\begin{split}
				\frac{\partial}{\partial t}\int_M\exp\frac{F^{*2}(D(u+tv))}{2}d\mu\mid_{t=0}=&\int_M\mathcal{V}(u)[Du(\nabla v)-C_{\nabla u}(\nabla u,\nabla u,\nabla v)]d\mu\\
				=&\int_M\mathcal{V}(u)du(\nabla v)d\mu.
			\end{split}
		\end{eqnarray}
		The conclusion holds due to the arbitrariness of $v$.
	\end{proof}

	As in the Riemannian case, the definition of an exponentially harmonic function $u$ on a Finsler metric measure spaces $(M,F,d\mu)$ is 
	\begin{eqnarray}
		\tilde\Delta_{\mu}u=0.
	\end{eqnarray}
	
	A direct computation shows that 
	\begin{eqnarray}
		\tilde\Delta_{\mu}u:=\mathrm{div}_{\mu}(\mathcal{V}(u)Du)=\mathcal{V}(u)\left(\Delta u+(Du)^2(\nabla^2u)\right).
	\end{eqnarray}
	By defining 
	\begin{eqnarray}
		\hat\Delta u=\frac{\tilde\Delta_{\mu}u}{\mathcal{V}(u)}=\Delta u+(du)^2(\nabla^2 u),
	\end{eqnarray}
	a function $u$ is exponentially harmonic on $(M,F,d\mu)$ if 
	\begin{eqnarray}
		\tilde\Delta_{\mu}u=\hat\Delta u=0,
	\end{eqnarray}
	or equivalently, 
	\begin{eqnarray}
		\Delta^{\nabla u}u+(Du)^2(\nabla^{\nabla u 2}u)=0.
	\end{eqnarray}
	Furthermore, we denote the second-order quasilinear elliptic operator $\Delta^{\nabla u}+(Du)^2\nabla^{\nabla u 2}$ by $\hat \Delta^{\nabla u}$.
	
	Let $e(u)=F^2(\nabla u)=F^{*2}(Du)=g^{ij}(\nabla u)u_iu_j$ be the energy density function. One can directly verify that
	\begin{eqnarray}
		\begin{split}
			\hat\Delta^{\nabla u}e=&\Delta^{\nabla u}e+(Du)^2(\nabla^{\nabla u 2}e)\\
			=&\mathrm{tr}_{\nabla u}\nabla^{\nabla u2}e-g_{\nabla u}(\nabla^{\nabla u} \tau,\nabla^{\nabla u} e)+(Du)^2\nabla^{\nabla u 2}e.
		\end{split}
	\end{eqnarray}
	Using the Ricci-type identity from \cite{Shen2018}, this simplifies to 
	\begin{eqnarray}\label{hatDe-1}
		\begin{split}
			\frac12\hat\Delta^{\nabla u}e=&Du(\nabla^{\nabla u}\mathrm{tr}_{\nabla u}\nabla^2u)+Ric(\nabla u)+\|\nabla^2u\|_{HS(\nabla u)}^2\\
			&-(Du\otimes D\tau)(\nabla^2 u)+(Du)^3(\nabla^3u)+F^2(Du(\nabla^2u)).
		\end{split}
	\end{eqnarray}
	Noticing that 
	\begin{eqnarray}
		-(Du\otimes D\tau)(\nabla^2 u)=-\dot S(\nabla u)+(Du)^2(\nabla^{\nabla u 2}\tau),
	\end{eqnarray}
	and using the definition of the weighted Ricci curvature $Ric^{\infty}$, we derive the following Bochner-type formula for exponentially harmonic functions from \eqref{hatDe-1}.
	\begin{eqnarray}\label{hatDe-2}
		\hat\Delta^{\nabla u}e=2Ric^{\infty}(\nabla u)+2\|\nabla^2u\|_{HS(\nabla u)}^2-\frac12F^{*2}(De(u)).
	\end{eqnarray}
	Moreover, if $\phi: \mathbb{R}\rightarrow\mathbb{R}$ is a $C^2$ function,  we obtain the equivalent expression for $\hat\Delta^{\nabla u}(\phi\cdot u)$ as
	\begin{eqnarray*}
		\begin{split}
			\hat\Delta(\phi\cdot u)=&\phi'\mathrm{tr}_{\nabla u}\nabla^2 u+\phi''F^2(\nabla u)-\phi'Du(\nabla^{\nabla u}\tau)+\phi'(Du)^2(\nabla^2u)+\phi''F^4(\nabla u)\\
			=&\phi'\hat\Delta u+\phi''(e(u)+e^2(u)),
		\end{split}
	\end{eqnarray*}
	so that
	\begin{eqnarray}\label{hatDphiu}
		\hat\Delta(\phi\cdot u)=\phi''\left(e(u)+e^2(u)\right),
	\end{eqnarray}
	provided $u$ is an exponentially harmonic function on $(M,F,du)$.

	\section{Proof of the main theorem}
	
	Let $x_0$ be a fixed point on $(M,F,d\mu)$ and let $B_a(x_0)$ denote the closed forward geodesic ball of radius $a$ centered at $x_0$. Denoting the forward distance from $x_0$ by $r$, we choose $\phi=4r^2$ and select $b$ sufficiently large so that $\sup\{\phi\cdot u(x) \mid x\in B_{2a}(x_0)\}<b^2<\infty$ on $B_{2a}(x_0)$. Consider the function
	$$H=\frac{(a^2-r^2)^2}{b^2-u^2}e(u),$$
	defined on the forward geodesic ball $B_a(x_0)$. Since $H$ vanishes on the boundary of $B_a(x_0)$, we assume that $H$ attains its maximum at an interior point $p$. Consequently, at $p$, the following conditions hold.
	\begin{eqnarray}\label{dlogH=0}
		0=d\log H=-2\frac{d(r^2)}{a^2-r^2}+\frac{de}{e}+\frac{d(u^2)}{b^2-u^2},
	\end{eqnarray}
	and
	\begin{eqnarray}\label{hatDlogH-1}
		0\geq& \hat\Delta^{\nabla u}(\log H).
	\end{eqnarray}
	\eqref{hatDlogH-1} is valid because the elliptic part of $\hat\Delta^{\nabla u}$ is a positive definite matrix.
	
	Since 
	\begin{eqnarray}
		\begin{split}
			d^2\log H=&-2\frac{d^2(r^2)}{a^2-r^2}-2\frac{d(r^2)\otimes d(r^2)}{(a^2-r^2)^2}+\frac{d^2 e}{e}\\
			&-\frac{de\otimes de}{e^2}+\frac{d^2(u^2)}{b^2-u^2}+\frac{d(u^2)\otimes d(u^2)}{(b^2-u^2)^2},
		\end{split}
	\end{eqnarray}
	the right-hand side (RHS) of \eqref{hatDlogH-1} becomes
	\begin{eqnarray*}
		\begin{split}
			\hat\Delta^{\nabla u}(\log H)=&-2\frac{\mathrm{tr}_{\nabla u}\nabla^{\nabla u2}(r^2)}{a^2-r^2}-2\frac{F_{\nabla u}^2(\nabla^{\nabla u}(r^2))}{(a^2-r^2)^2}+\frac{\mathrm{tr}_{\nabla u}(\nabla^{\nabla u2}e)}{e}\\
			&-\frac{F_{\nabla u}^2(\nabla^{\nabla u}e)}{e^2}+\frac{\mathrm{tr}_{\nabla u}\nabla^2(u^2)}{b^2-u^2}+\frac{F_{\nabla u}^2(\nabla^{\nabla u}(u^2))}{(b^2-u^2)^2}\\
			&+2\frac{g_{\nabla u}(\nabla^{\nabla u}\tau,\nabla^{\nabla u}(r^2))}{a^2-r^2}-\frac{g_{\nabla u}(\nabla^{\nabla u}\tau,\nabla^{\nabla u}e)}{e^2}-\frac{D(u^2)(\nabla^{\nabla u}\tau)}{b^2-u^2}\\
			&-2\frac{(Du)^2(\nabla^{\nabla u2}(r^2))}{a^2-r^2}-2\frac{[Du(\nabla^{\nabla u}(r^2))]^2}{(a^2-r^2)^2}+\frac{(Du)^2(\nabla^{\nabla u2}e)}{e}\\
			&-\frac{[Du(\nabla^{\nabla u}e)]^2}{e^2}+\frac{(Du)^2(\nabla^2(u^2))}{b^2-u^2}+\frac{[Du(\nabla(u^2))]^2}{(b^2-u^2)^2},
		\end{split}
	\end{eqnarray*}
	which is equivalent to
	\begin{eqnarray}\label{hatDlogH-2}
		\begin{split}
			\hat\Delta^{\nabla u}(\log H)=&\frac{-2\hat\Delta^{\nabla u}r^2}{a^2-r^2}-\frac{2F^{*2}(\tilde d r^2)}{(a^2-r^2)^2}+\frac{\hat\Delta^{\nabla u}e}{e}\\
			&-\frac{F^{*2}(\tilde de)}{e^2}+\frac{\hat\Delta u^2}{b^2-u^2}+\frac{F^{*2}(\tilde du^2)}{(b^2-u^2)^2},
		\end{split}
	\end{eqnarray}
	where the notation $\tilde d$ is defined by $F^{*2}(\tilde df)=F^{*2}(Df)+[Df(\nabla u)]^2$.
	
	From \eqref{hatDphiu}, we obtain that
	\begin{eqnarray}\label{hatDu2}
		\hat\Delta u^2\geq 8e(u)(1+e(u)).
	\end{eqnarray}
	The curvature condition and the Cauchy-Schwartz inequality $4\|\nabla^2 u\|_{HS(\nabla u)}^2e(u)\geq F^{*2}(De(u))$ allow us to refine \eqref{hatDe-2} to 
	\begin{eqnarray}\label{hatDe}
		\hat\Delta^{\nabla u}e\geq2\|\nabla^2u\|_{HS(\nabla u)}^2-\frac12F^{*2}(De(u))\geq\frac{F^{*2}(De)}{2e}(1-e).
	\end{eqnarray}
	From \eqref{dlogH=0}, it follows that 
	\begin{eqnarray}\label{Ftildede}
		\begin{split}
			\frac{F^{*2}(\tilde de)}{e^2}=&\frac{1}{e^2}\{F^{*2}(De)+[De(\nabla u)]^2\}=F^{*2}(\frac{De}{e})+\left[\frac{De}{e}(\nabla u)\right]^2\\
			=&\frac{4F^{*2}(\tilde dr^2)}{(a^2-r^2)^2}+\frac{F^{*2}(\tilde du^2)}{(b^2-u^2)^2}-\frac{4(1+e)g_{\nabla u}(\nabla u^2,\nabla^{\nabla u}r^2)}{(a^2-r^2)(b^2-u^2)}.		
		\end{split}
	\end{eqnarray}
	Thus, combining \eqref{hatDlogH-1}, \eqref{hatDlogH-2}-\eqref{Ftildede}, we obtain
	\begin{eqnarray}\label{hatDlogH-3}
		\begin{split}
			0\geq\hat\Delta^{\nabla u}(\log H)=&\frac{-2\hat\Delta^{\nabla u}r^2}{a^2-r^2}-\frac{6F^{*2}(\tilde d r^2)}{(a^2-r^2)^2}+\frac{F^{*2}(De)}{2e^2}(1-e)\\
			&-\frac{4(1+e)g_{\nabla u}(\nabla u^2,\nabla^{\nabla u}r^2)}{(a^2-r^2)(b^2-u^2)}+\frac{8e(u)(1+e(u))}{b^2-u^2}.
		\end{split}
	\end{eqnarray}
	
	It follows from Theorem \ref{thm-LapComp-1} that
	\begin{eqnarray}\label{hatDr2}
		\begin{split}
			\hat\Delta^{\nabla u}r^2=&2r\Delta^{\nabla u}r+2F^{*2}(Dr)+2\left[Du(\nabla^{\nabla u}r)\right]^2+2r(Du)^2(\nabla^{\nabla u2}r)\\
			\leq&2\left[r(\Delta^{\nabla u}r)+\alpha\right]+2e(u)\left[\alpha+r(|\Delta^{\nabla u}r|+|g_{\nabla u}(\nabla^{\nabla u}\tau,\nabla^{\nabla u}r)|)\right]\\
			\leq&2\left[C(N,\alpha)+\alpha+\sqrt{\alpha}K_0r\right]+2e(u)\left[\alpha+C(N,\alpha)+2\alpha K_0+\sqrt{\alpha}K_0r\right]\\
			\leq&C_1(1+C_2r)(1+e(u)),
		\end{split}
	\end{eqnarray}
	where $C_1$ and $C_2$ are two constants depending on the dimension $n$, the misalignment $\alpha$, and the bound of non-Riemannian curvatures $K_0$.
	Moreover, \eqref{dlogH=0} implies that
	\begin{eqnarray}\label{Fde}
		\begin{split}
			\frac{F^{*2}(De)}{e^2}=&\frac{16\alpha r^2}{(a^2-r^2)^2}-\frac{16rudr(\nabla u)}{(a^2-r^2)(b^2-u^2)}+\frac{4u^2e}{(b^2-u^2)^2}\\
			\leq&\frac{32\alpha r^2}{(a^2-r^2)^2}+\frac{8u^2e}{(b^2-u^2)^2}.
		\end{split}
	\end{eqnarray}
	
	Additionally, we have
	\begin{eqnarray}\label{Fdtilder2}
		\begin{split}
			F^{*2}(\tilde d r^2)=&F^2_{\nabla u}(\nabla^{\nabla u}r^2)+\left[Dr^2(\nabla u)\right]^2\\
			\leq& 4r^2F^{*2}(Dr)+4r^2F^{*2}(Dr)F^2(\nabla u)=4\alpha r^2(1+e(u)),
		\end{split}
	\end{eqnarray}
	and
	\begin{eqnarray}\label{gdu2dr2}
		\begin{split}
			\frac{4(1+e)g_{\nabla u}(\nabla u^2,\nabla^{\nabla u}r^2)}{(a^2-r^2)(b^2-u^2)}\leq&\frac{16(1+e)\sqrt{\alpha}urF(\nabla u)}{(a^2-r^2)(b^2-u^2)}\\
			\leq& \frac{16\alpha(1+e)r^2u^2}{(a^2-r^2)^2(b^2-u^2)}+\frac{4e(u)(1+e(u))}{b^2-u^2}.
		\end{split}
	\end{eqnarray}
	
	Substituting \eqref{hatDr2}-\eqref{gdu2dr2} into \eqref{hatDlogH-3}, we conclude that
	\begin{eqnarray}\label{hatDlogH-4}
		\begin{split}
			0\geq&-\frac{2C_1(1+C_2r)(1+e(u))}{a^2-r^2}-\frac{24\alpha r^2(1+e)}{(a^2-r^2)^2}-\frac{16\alpha r^2e}{(a^2-r^2)^2}\\
			&-\frac{4u^2e^2}{(b^2-u^2)^2}-\frac{16\alpha(1+e)r^2u^2}{(a^2-r^2)^2(b^2-u^2)}+\frac{4e(u)(1+e(u))}{b^2-u^2}.
		\end{split}
	\end{eqnarray}
	By choosing $b$ sufficiently large, we ensure that 
	\begin{eqnarray}
		\frac{4u^2}{(b^2-u^2)^2}<\frac{2}{(b^2-u^2)},
	\end{eqnarray}
	which, when combined with \eqref{hatDlogH-4}, yields
	\begin{eqnarray}
		\begin{split}
			0\geq\frac{2e^2}{b^2-u^2}&-e\left(\frac{2C_1(1+C_2r)}{a^2-r^2}+\frac{40\alpha r^2}{(a^2-r^2)^2}+\frac{16\alpha r^2u^2}{(a^2-r^2)^2(b^2-u^2)}\right)\\
			&-\left(\frac{2C_1(1+C_2r)}{a^2-r^2}+\frac{24\alpha r^2}{(a^2-r^2)^2}+\frac{16\alpha r^2u^2}{(a^2-r^2)^2(b^2-u^2)}\right)\\ 
			=:&2(Ae^2-Be-E).
		\end{split}
	\end{eqnarray}
	Therefore, it follows that
	\begin{eqnarray}
		e(u)(p)\leq\max\left\{\frac{2B}{A},\frac{2\sqrt{E}}{A}\right\}.
	\end{eqnarray}
	That is,
	\begin{eqnarray}
		\frac{a^2-r_0^2}{b^2-u^2(p)}e(u)(p)\leq\max\left\{C_4(a^3+a^2+1),C_5(a^{3+\frac12}+a^2+1)\right\},
	\end{eqnarray}
	where $r_0=d(x_0,p)$ is the forward distance from $x_0$ to $p$, and $C_4$, $C_5$ are constants depending on $b$, $n$, $\alpha$, $K_0$ and $\sup u^2$.
	
	Noticing that $H(x)\leq H(p)$ for all $x\in B_a(x_0)$, we easily deduce that
	\begin{eqnarray}\label{e<Ca}
		\frac{(a^2-r^2)^2e(u)}{b^2-u^2}\leq C\left(a^{3+\frac12}+a^2+1\right),
	\end{eqnarray}
	for all $x\in B_a(x_0)$. Then for any fixed point $x$, taking the limit as $a\rightarrow \infty$ in \eqref{e<Ca} shows that
	\begin{eqnarray}
		e(u)=0,
	\end{eqnarray}
	for all $x$ on $M$. Thus, $u$ must be a constant.


\section*{Acknowledgments}

The author is grateful to the reviewers for their careful reviews and valuable comments. 

\section*{Declarations}

\subsection*{Ethical Approval}
Ethical Approval is not applicable to this article as no human or animal studies in this study.

\subsection*{Funding} 
The first author is supported partially by the NNSFC (Nos. 12001099, 12271093, 12271093).

\subsection*{Data availability statement}
Data sharing is not applicable to this article as no new data were created or analyzed in this study.

\subsection*{Materials availability statement}
Materials sharing is not applicable to this article.

\subsection*{Conflict of interest/Competing interests}
All authors disclosed no relevant relationships.

\end{document}